\documentclass[12pt,a4paper]{article}
\usepackage{amssymb,amsmath,amsthm}

\usepackage{url,color}
\usepackage{times,lscape}
\def\de{{\rm d}}

\newtheorem{theorem}{Theorem}[section]
\newtheorem{lemma}{Lemma}[section]

\newtheorem{remark}{Remark}[section]
\newtheorem{assumption}{Assumption}[section]
\numberwithin{equation}{section}

\title{Divergences Test Statistics for Discretely Observed Diffusion Processes}

\author{Alessandro De Gregorio\\ \small Department of Statistics, Probability and Applied Statistics\\ \small P.le Aldo Moro 5, 00185 Rome-Italy\\ \small alessandro.degregorio@uniroma1.it\and Stefano M. Iacus\\ \small Department of Economics, Business and Statistics\\ \small Via Conservatorio 7, 20124 Mlan-Italy\\ \small stefano.iacus@unimi.it}

\begin{document}
\maketitle

\begin{abstract}
In this paper we propose the use of $\phi$-divergences as test
statistics to verify  simple hypotheses about a one-dimensional
parametric diffusion process $\de X_t = b(X_t, \theta)\de t +
\sigma(X_t, \theta)\de W_t$, from discrete observations
$\{X_{t_i}, i=0, \ldots, n\}$ with $t_i = i\Delta_n$, $i=0, 1,
\ldots, n$, under the asymptotic scheme $\Delta_n\to0$,
$n\Delta_n\to\infty$ and $n\Delta_n^2\to 0$. The class of
$\phi$-divergences is wide and includes several special members
like Kullback-Leibler, R\'enyi, power and $\alpha$-divergences. We
derive the asymptotic distribution of the test statistics based on
$\phi$-divergences. The limiting law takes different forms
depending on the regularity of $\phi$. These convergence differ
from the classical results for independent and identically
distributed random variables. Numerical analysis is used to show
the small sample properties of the test statistics in terms of
estimated level and power of the test.
\end{abstract}
\noindent {\bf keywords:} diffusion processes, empirical level,
hypotheses testing, $\phi$-divergences, $\alpha$-divergences

\section{Introduction}
We consider the problem of parametric testing using $\phi$-divergences.
Let $X$ be a r.v. and $f(X,\theta)$ and $g(X,\theta)$, $\theta\in
\Theta$ two families of probability densities on the same
measurable space. The $\phi$-divergences are defined as
$D_\phi(f,g)=E_\theta \phi\left(f(X)/g(X)\right)$,
where $E_\theta$ is the expected value with respect to $P_\theta$, the true law of
the observations. Because we focus the attention on the
use of divergences for hypotheses testing, we will use a
simplified notation: let $\theta$ and $\theta_0$ two points in the
interior of $\Theta$ and define the divergence as
\begin{equation}
\label{eq:div}
D_\phi(\theta,\theta_0)=E_{\theta_0}\phi\left(\frac{p(X,\theta)}{p(X,\theta_0)}\right)
\end{equation}
In equation \eqref{eq:div} the density $\{ p(X,\theta),
\theta\in\Theta\}$ is a same family of probability densities and
$\phi(\cdot)$ is a function with the minimal property that
$\phi(1)=0$. Examples of divergences of the form
$D_\alpha(\theta,\theta_0)=D_{\phi_\alpha}(\theta,\theta_0)$ are
the $\alpha$-divergences, defined by means of the following
function
$$\phi_\alpha(x) =\frac{4(1-x^\frac{1+\alpha}{2})}{1-\alpha^2}, \quad  -1 < \alpha <1$$
Note that $D_\alpha(\theta_0, \theta) = D_{-\alpha}(\theta,
\theta_0)$. The class of $\alpha$-divergences has been widely
studied in statistics (see, e.g.,  Csisz\'ar, 1967 and Amari,
1985) and it is a family of divergences which includes several
members of particular interest. For example, in the limit as
$\alpha\to -1$, $D_{-1}(\theta, \theta_0)$ reduces to the
well-known
 Kullback-Leibler measure
$$
D_{-1}(\theta, \theta_0) = - E_{\theta_0} \log\left(\frac{p(X,\theta)}{p(X,\theta_0)}\right)
$$
while as $\alpha\to0$,  the Hellinger distance (see, e.g., Beran,
1977, Simpson, 1989) emerges
$$
D_0(\theta, \theta_0) = \frac12 E \left(\sqrt{p(X,\theta)} - \sqrt{p(X,\theta_0)}\right)^2
$$
As noticed in Chandra and Taniguchi (2006), the $\alpha$-divergence is also equivalent to
the R\'enyi's divergence (R\'enyi, 1961) defined, for $\alpha\in(0,1)$, as
$$
R_\alpha(\theta, \theta_0)= \frac{1}{1-\alpha} \log E_{\theta_0} \left( \frac{p(X,\theta)}{p(X,\theta_0)} \right)^\alpha
$$
from which is easy to see that in the limit as $\alpha\to1$,
$R_\alpha$ reduces to the Kullback-Leibler divergence. The
transformation $\psi(R_\alpha) = (\exp\{(\alpha-1)R_\alpha
-1\}/(1-\alpha)$ returns the power-divergence studied in Cressie
and Read (1984). Liese and Vajda (1987) provide extensive study of
a modified version of $R_\alpha$ and Morales {\it et al.} (1997)
consider divergences with convex $\phi(\cdot)$  for independent
and identically distributed (i.i.d) observations; for example the
power-divergences $D_{\phi_\lambda}(\theta,\theta_0)$ with
\begin{equation}
\phi_\lambda(x)=\frac{x^\lambda-\lambda(x-1)-1}{\lambda(\lambda-1)}
\end{equation}
and $\lambda\in \mathbb{R}-\{0,1\}$.

In this paper we focus our attention on the $\phi$-divergences
$D_\phi(\theta,\theta_0)$, defined as in \eqref{eq:div}, for
one-dimensional diffusion process $\{X_t, t\in [0,T]\}$, solution
of the following stochastic differential equation
\begin{equation}
dX_t=b(\alpha,X_t)dt+\sigma(\beta,X_t)dW_t,\quad X_0=x_0,
\label{eq:diff}
\end{equation}
where $W_t$ is a Brownian motion,
 $\theta=(\alpha,\beta)\in \Theta_\alpha\times\Theta_\beta=\Theta$,
 where $\Theta_\alpha$ and $\Theta_\beta$ are
 respectively compact convex subset of $\mathbb{R}^p$ and $\mathbb{R}^q$.
  We assume that the process $X_t$ is ergodic for every $\theta$ with invariant
 law $\mu_\theta$. Furthermore $X_t$ is observed at discrete
 times $t_i=i\Delta_n,i=0,1,2,...,n,$ where $\Delta_n$ is the length of the steps.
 We indicate the observations with $\bold{X}_n=\{X_{t_i}\}_{0\leqslant i \leqslant n}$.
 The asymptotic is $\Delta_n\to 0, n\Delta_n\to \infty$ and $n\Delta_n^2\to 0$ as $n\to \infty$.

We study the properties of  the estimated $\phi$-divergence
$\mathbb D_\phi(\tilde\theta_n(\bold{X}_n),\theta_0)$, for
discretely observed diffusion processes, defined as
$$
\mathbb{D}_\phi(\tilde\theta_n(\bold{X}_n),\theta_0)=\phi\left(\frac{f_n({\bf
X}_n,\tilde\theta_n({\bf X}_n))}{f_n({\bf X}_n,\theta_0)}\right)
$$
where $f_n(\cdot, \cdot)$  is the approximated likelihood proposed
by Dacunha-Castelle and Florens-Zmirou (1986) and
$\tilde\theta_n(\bold{X}_n)$ is any consistent, asymptotically
normal and efficient estimator of $\theta$. We prove that, for
$\phi(\cdot)$ functions which satisfying three different
regularity conditions,  the statistic $\mathbb{D}_\phi$ converge
weakly to three different functions of the  $\chi_{p+q}^2$  random
variable. This result differs from the case of i.i.d. setting.

Up to our knowledge the only result concerning the use of
divergences for discretely observed diffusion process is due to
Rivas {\it et al.} (2005) where they consider the model of
Brownian motion with drift $\de X_t = a \de t + b \de W_t$ where
$a$ and $b$ are two scalars. In that case, the exact likelihood of
the observations is available in explicit form and is the gaussian
law. Conversely, in the general setup of this paper, the
likelihood of the process in \eqref{eq:diff} is known only for
three particular stochastic differential equations, namely the
Ornstein-Uhlembeck diffusion, the geometric Brownian motion and
the Cox-Ingersoll-Ross model. In all other cases, the likelihood
has to be approximated. We choose the approximation due to
Dacunha-Castelle and Florens-Zmirou (1986) and, to derive a proper
estimator, we use the local gaussian approximation proposed by
Yoshida (1992) although our result holds for any consistent and
asymptotically Gaussian estimator. This approach has been
suggested by the work on Akaike Information Criteria by Uchida and
Yoshida (2005).

For continuous time observations from diffusion processes, Vajda (1990) considered the model $\de X(t) = - b(t)X_t \de t + \sigma(t) \de W_t$;  K\"uchler and S\o rensen (1997) and Morales {\it et al.} (2004) contain several results on the likelihood ratio test statistics and R\'enyi statistics for exponential family of diffusions. Explicit derivations of the R\'enyi information on the invariant law of  ergodic diffusion processes have been presented in De Gregorio and Iacus (2007).
For small diffusion processes, with continuous time observations, information criteria have been derived in Uchida and Yoshida (2004) using Malliavin calculus.

The problem of testing statistical hypotheses from general diffusion processes is still a developing stream of research.
Kutoyants (2004) and Dachian and Kutoyants (2008) consider the problem of testing statistical hypotheses for ergodic diffusion models in continuous time; Kutoyants (1984) and  Iacus and Kutoyants (2001) consider parametric and semiparametric  hypotheses testing for small diffusion processes; Negri and Nishiyama (2007a, b) propose a non parametric test based on score marked empirical process for both continuous and discrete time observation from small diffusion processes further extended to the ergodic case in Masuda {\it et al.} (2008). Lee and Wee (2008) considered the parametric version of the same test statistics for a simplified model.

A\"it-Sahalia (1996, 2008), Giet and Lubrano (2008) and Chen {\it
et al.} (2008) proposed tests based on the several distances
between parametric and nonparametric estimation of the invariant
density of discretely observed ergodic diffusion processes. The
present paper complements the above references.

The paper is organized as follows. Section \ref{sec:notation}
introduces notation and regularity assumptions. Section
\ref{sec:results} states the main result. Section
\ref{sec:applications} contains numerical experiments to test the
small sample performance of the proposed test statistics in terms
of empirical level and empirical power under some alternatives.
The proofs are contained in  Section \ref{sec:proofs}.

\section{Assumptions on diffusion model}\label{sec:notation}
We consider the family of one-dimensional diffusion processes $\{X_t, t\in [0,T]\}$, solution to
\begin{equation}
dX_t=b(\alpha,X_t)dt+\sigma(\beta,X_t)dW_t,\quad X_0=x_0,
\label{eq:sde1}
\end{equation}
where $W_t$ is a Brownian motion.
Let $\theta=(\alpha,\beta)\in \Theta_\alpha\times\Theta_\beta=\Theta$, where $\Theta_\alpha$ and $\Theta_\beta$ are
 respectively compact convex subset of $\mathbb{R}^p$ and $\mathbb{R}^q$. Furthermore we assume that the drift function $b:\mathbb{R\times }\Theta_\alpha\to \mathbb{R}$ and the diffusion coefficient $\sigma:\mathbb{R }\times\Theta_\beta\to \mathbb{R}$ are known apart from the parameters $\alpha$ and $\beta$.
 We assume that the process $X_t$ is ergodic for every $\theta$ with invariant
 law $\mu_\theta$. The process $X_t$ is observed at discrete
 times $t_i=i\Delta_n, i=0,1,2,...,n,$ where $\Delta_n$ is the length of the steps.
 We indicate the observations with $\bold{X}_n=\{X_{t_i}\}_{0\leqslant i \leqslant n}$.
 The asymptotic is $\Delta_n\to 0, n\Delta_n\to \infty$ and $n\Delta_n^2\to 0$ as $n\to \infty$.

In the definition of the $\phi$-divergence \eqref{eq:div} the likelihood of the process is
need, but as noted in the Introduction, this is usually not know. There are several ways to
approximate the likelihood of a discretely observed diffusion process
(for a review see, e.g., Chap. 3, Iacus, 2008). In this paper, we use the approximation
 proposed by Dacunha-Castelle and Florens-Zmirou (1986) although our result hold true
 (with some adaptations of the proofs) for other approximations, like, e.g. the one based
 on Hermite polynomial expansion by A\"it-Sahalia (2002). To write it in explicit way, we use the same setup as in Uchida and Yoshida (2005). We introduce the following functions
$$s(x,\beta)=\int_0^x\frac{du}{\sigma(\beta,u)},\quad
B(x,\theta)=\frac{b(\alpha,x)}{\sigma(\beta,x)}-\frac{\sigma'(\beta,x)}{2}$$
$$\widetilde{B}(x,\theta)=B(s^{-1}(\beta,x),\theta), \quad
\widetilde{h}(x,\theta)=\widetilde{B}^2(x,\theta)+\widetilde{B}'(x,\theta)$$
The following set of assumptions ensure the good behaviour of the approximated likelihood
and the existence of a weak solution of \eqref{eq:sde1}
\begin{assumption}\label{ass:1}
[Regularity on the process]
\begin{itemize}
\item[i)] There exists a constant $C$ such that
$$|b(\alpha_0,x)-b(\alpha_0,y)|+|\sigma(\beta_0,x)-\sigma(\beta_0,y)|\leq C|x-y|.$$

\item[ii)] $\inf_{\beta,x}\sigma^2(\beta,x)>0$.

\item[iii)] The process $X$ is ergodic for every $\theta$ with
invariant probability measure $\mu_\theta$. All polynomial moments
of $\mu_\theta$ are finite.

\item[iv)] For all $m\geq 0$ and for all $\theta$, $\sup_t
E|X_t|^m<\infty$.

\item[v)] For every $\theta$, the coefficients $b(\alpha,x)$ and
$\sigma(\beta,x)$ are twice differentiable with respect to $x$ and
the derivatives are polynomial growth in $x$, uniformly in
$\theta$.

\item[vi)] The coefficients $b(\alpha,x)$ and $\sigma(\beta,x)$
and all their partial derivatives respect to $x$ up to order 2 are
three times differentiable respect to $\theta$ for all $x$ in the
state space. All derivatives respect to $\theta$ are polynomial
growth in $x$, uniformly in $\theta$.
\end{itemize}
\end{assumption}

\begin{assumption}\label{ass:2}
[Regularity for the approximation]
\begin{itemize}
\item[i)] $\widetilde{h}(x,\theta)=O(|x|^2)$ as $x\to \infty$.

\item[ii)] $\inf_x\widetilde{h}(x,\theta)>-\infty$ for all
$\theta$.

\item[iii)] $\sup_\theta\sup_x|\widetilde{h}^3(x,\theta)|\leq M
<\infty$.

\item[iv)] There exists $\gamma>0$ such that for every $\theta$
and $j=1,2,$
$|\widetilde{B}^j(x,\theta)|=O(|\widetilde{B}(x,\theta)|^\gamma)$
as $|x|\to\infty$.
\end{itemize}
\end{assumption}

\begin{assumption}\label{ass:3}
[Identifiability]
The coefficients $b(\alpha,x)=b(\alpha_0,x)$ and
$\sigma(\beta,x)=\sigma(\beta_0,x)$ for $\mu_{\theta_0}$ a.s. all
x then $\alpha=\alpha_0$ and $\beta=\beta_0$.
\end{assumption}
Under Assumptions \ref{ass:1} and \ref{ass:2} Dacunha-Castelle and
Florens-Zmirou (1986) introduced the following approximation of
transition density $f$ of the process $X$ from $y$ to $x$ at lag
$t$
\begin{equation}
f(x,y,t,\theta)=\frac{1}{\sqrt{2\pi t}\sigma(y,\beta)}\exp\left\{-\frac{S^2(x,y,\beta)}{2t}+H(x,y,\theta)+t\tilde g(x,y,\theta)\right\}
\label{eq:approx1}
\end{equation}
and its logarithm
\begin{eqnarray*}
l(x,y,t,\theta)&=&-\frac12\log(2\pi t)-\log\sigma(y,\beta)-\frac{S^2(x,y,\beta)}{2t}+H(x,y,\theta)+t\tilde g(x,y,\theta)
\end{eqnarray*}
where
$$S(x,y,\beta)=\int_x^y\frac{du}{\sigma(u,\beta)}$$
$$H(x,y,\theta)=\int_x^y \left\{\frac{b(\alpha,u)}{\sigma^2(\beta,u)}-\frac12\frac{\sigma'(\beta,u)}{\sigma(\beta,u)}\right\}du$$
$$\tilde g(x,y,\theta)=-\frac12\left\{C(x,\theta)+C(y,\theta)+\frac13B(x,\theta)B(y,\theta)\right\}$$
$$C(x,\theta)=\frac13B^2(x,\theta)+\frac12B'(x,\theta)\sigma(x,\beta)$$
The approximated likelihood and log-likelihood functions of the
observations  $\bold{X}_n$ become respectively
$$f_n(\bold{X}_n,\theta)=\prod_{i=1}^n f(\Delta_n,X_{t_{i-1}},X_{t_{i}},\theta)$$
$$l_n(\bold{X}_n,\theta)=\sum_{i=1}^n l(\Delta_n,X_{t_{i-1}},X_{t_{i}},\theta)$$

\section{Construction of the test statistics and results}\label{sec:results}
Consider the divergence defined in \eqref{eq:div} and let
$\phi(\cdot)$ be such that $\phi(1)=0$ and,
 when they exist, define  $C_\phi = \phi'(1)$ and $K_\phi=\phi''(1)$.
We consider three different setup
\begin{assumption}\label{ass:C1}
  $C_\phi\neq 0$ is a finite constant depending only on $\phi$ and independent of $\theta$;
\end{assumption}
\begin{assumption}\label{ass:C2}
 $C_\phi=0$ and $K_\phi\neq 0$ is a finite constant
depending only on $\phi$ and independent of $\theta$;
\end{assumption}
\begin{assumption}\label{ass:C3}
 $C_\phi\neq 0$ and
$K_\phi\neq 0$ are finite constants depending only on $\phi$ and
independent of $\theta$;
\end{assumption}

\begin{remark}The above Assumptions are not so strong. In fact, for example the $\alpha$-divergences
 $D_{\phi_\alpha}(\theta,\theta_0)$ satisfy the Assumptions \ref{ass:C1}
and \ref{ass:C3}, while for the power-divergences
$D_{\phi_\lambda}(\theta,\theta_0)$ it's easy to verify that
$C_\phi=\phi'(1)=0$.
\end{remark}
 Clearly, the quantity
$D_\phi(\theta,\theta_0)$ measures the discrepancy between
$\theta$ and the true value of the parameter $\theta_0$ and is an
ideal candidate to construct a test statistics. Let
$\tilde\theta_n(\bold{X}_n)$ be any consistent estimator of
$\theta_0$ and such that
\begin{equation}
\label{eq:conv0}
\Gamma^{-1/2}(\tilde\theta_n(\bold{X}_n)-\theta_0)\stackrel{d}{\to}N(0,\mathcal
I(\theta_0)^{-1})
\end{equation}
where $\mathcal I(\theta_0)$ is the positive definite and invertible
Fisher information matrix at $\theta_0$  equal to
$$\mathcal I(\theta_0)=\left(%
\begin{array}{cc}
  (\mathcal I_b^{kj}(\theta_0))_{k,j=1,...,p} & 0 \\
  0 & (\mathcal I_\sigma^{kj}(\theta_0))_{k,j=1,...,q} \\
\end{array}%
\right)$$ where
$$\mathcal I_b^{kj}(\theta_0)=\int\frac{1}{\sigma^2(\beta_0,x)}\frac{\partial b(\alpha_0,x)}{\partial\alpha_k}\frac{\partial b(\alpha_0,x)}{\partial\alpha_j}
\mu_{\theta_0}(dx)$$
$$\mathcal I_\sigma^{kj}(\theta_0)=2\int\frac{1}{\sigma^2(\beta_0,x)}\frac{\partial \sigma(\beta_0,x)}{\partial\beta_k}\frac{\partial \sigma(\beta_0,x)}{\partial\beta_j}
\mu_{\theta_0}(dx)$$ We indicate with $\Gamma$ the $(p+q)\times
(p+q)$ matrix
$$\Gamma=\left(%
\begin{array}{cc}
  \frac{1}{n\Delta_n}I_p & 0 \\
  0 &  \frac{1}{n}I_q\\
\end{array}%
\right)$$ and $I_p$ is the $p\times p$ identity matrix.
Using the approximated likelihood $f_n({\bf X}_n,\theta)$ and $f_n({\bf X}_n,\theta_0)$, the $\phi$-divergence in \eqref{eq:div} becomes
\begin{equation}
\label{eq:dvi2}
D_\phi(\theta,\theta_0)=E_{\theta_0} \phi\left(\frac{f_n({\bf X}_n,\theta)}{f_n({\bf X}_n,\theta_0)}\right)
\end{equation}
To construct a test statistics we replace $\theta$ by the  estimator $\tilde\theta_n(\bold{X}_n)$ and, having only one single observation of ${\bf X}_n$, i.e. only one observed trajectory, we estimate  \eqref{eq:dvi2} with
\begin{equation}
\label{eq:estdiv}
\mathbb{D}_\phi(\tilde\theta_n(\bold{X}_n),\theta_0)=\phi\left(\frac{f_n({\bf X}_n,\tilde\theta_n({\bf X}_n))}{f_n({\bf X}_n,\theta_0)}\right)
\end{equation}
Please notice that, conversely to the i.i.d. case, there is no integral in the definition of \eqref{eq:estdiv}. We will discuss this point after the presentation of the Theorem \ref{main}.
The proposed test for testing $H_0 : \theta = \theta_0$ versus $H_1 : \theta \neq \theta_0$ is realized as
$\mathbb{D}_\phi(\tilde\theta_n(\bold{X}_n),\theta_0) = 0$ versus $\mathbb{D}_\phi(\tilde\theta_n(\bold{X}_n),\theta_0)\neq 0$.
\begin{theorem}\label{main}
Under $H_0:\theta=\theta_0$, Assumptions \ref{ass:1}-\ref{ass:3}, convergence \eqref{eq:conv0}, we have that
\begin{itemize}
\item[i)] if  function $\phi(\cdot)$ satisfies Assumption \ref{ass:C1}, then
\begin{equation}\label{chi}
\mathbb{D}_\phi(\tilde\theta_n(\bold{X}_n),\theta_0)\stackrel{d}{\to}
C_\phi \chi^2_{p+q}
\end{equation}
\item[ii)] if  function $\phi(\cdot)$ satisfies Assumption \ref{ass:C2}, then
\begin{equation}\label{chi2}
\mathbb{D}_\phi(\tilde\theta_n(\bold{X}_n),\theta_0)\stackrel{d}{\to}
\frac{K_\phi}{2} Z_{p+q}
\end{equation}
where $\sqrt{Z_{p+q}}=\chi^2_{p+q}$.
 \item[iii)] if  function $\phi(\cdot)$ satisfies Assumption \ref{ass:C3}, then
\begin{equation}\label{chi3}
\mathbb{D}_\phi(\tilde\theta_n(\bold{X}_n),\theta_0)\stackrel{d}{\to}
\frac12 (C_\phi\chi^2_{p+q}+(C_\phi+K_\phi)Z_{p+q})
\end{equation}
\end{itemize}
\end{theorem}
\begin{remark}
It's clear that for $C_\phi=0$ from \eqref{chi3} we immediately
reobtain the convergence result \eqref{chi2}.
\end{remark}

\begin{remark}
If we consider the limits as $\alpha\to-1$ for $\phi_\alpha(x)$ of the $\alpha$-divergences, i.e. we consider the Kullback-Leibler divergence, we have
$$\phi(x) = \lim_{\alpha\to-1} \phi_\alpha(x) = -\log(x)$$
for which $C_\phi = -1$ and $K_\phi=1$. In that case, \eqref{chi3} reduces to the standard result for the likelihood ratio test statistics.

\end{remark}

The convergence in Theorem \ref{main} may appear somewhat strange if one thinks about the usual results on $\phi$-divergences for i.i.d. observations. The main difference in diffusion models, is that our estimate of the divergence has not the usual form of an expected value, i.e. it estimates the expected value with one observation only. This is why, in the i.i.d case,  the first term in the Taylor expansion of $D_\phi$ vanishes being the expected value of the score function,  while in our case it remains only the score function which, as usual, converges to a Gaussian random variable. For the same reason, in the second term of the Taylor expansion, in the i.i.d. case appears the expected value of the second order derivative which converges to the Fisher information and, in our case, we have not the expected value, hence the convergence to the square of the $\chi^2$ emerges.

If one wants to emulate the standard results for the i.i.d. case, it is still possible to work on the invariant density of the diffusion process. In that case, the $\phi$-divergence takes the usual form of the i.i.d. case because the invariant density have the explicit form. Indeed, let
$$ s(x,\theta) = \exp\left\{-2\int_{\tilde x}^x\frac{b(y,\theta)}{\sigma^2(y,\theta)}\de y\right\}, \quad
m(x,\theta) = \frac{1}{\sigma^2(x,\theta)s(x,\theta)}
$$
be the scale and speed functions of the diffusion, with $\tilde x$
some value in the state space of the diffusion process. Let $M =
\int m(x,\theta) \de x$, then $\pi(x,\theta) = m(x,\theta)/M$ is
the invariant density of the diffusion process. In this case, it
is possible to define the $\phi$-divergence as
$$
D_\phi(\tilde\theta_n, \theta_0) = \int
\phi\left(\frac{\pi(x,\tilde\theta_n)}{\pi(x,\theta_0)}\right)\pi(x,\theta_0)\de
x
$$
and the standard results follows.

\begin{remark}
In our application, to derive and estimator, we consider further the local gaussian approximation of the same transition density (see, Yoshida, 1992)
 \begin{equation}
 \label{eq:locgauss}
 g_n(\bold{X}_n,\theta)=\sum_{i=1}^ng_n(\Delta_n,X_{t_{i-1}},X_{t_{i}},\theta)
 \end{equation}
where
$$g(t,x,y,\theta)=-\frac12\log(2\pi t)-\log\sigma(\beta,x)-\frac{[y-x-tb(\alpha,x)]^2}{2t\sigma^2(\beta,x)}$$
The approximate maximum likelihood estimator $\hat\theta_n(\bold{X}_n)$ based on \eqref{eq:locgauss} is then defined as
\begin{equation}
\hat\theta_n(\bold{X}_n)=\arg\sup_\theta g_n(\bold{X}_n,\theta)
\label{eq:mle}
\end{equation}
Under the condition $n\Delta_n^2\to 0$ (see Theorem 1 in
Kessler, 1997) the estimator $\hat\theta_n(\bold{X}_n)$ in \eqref{eq:mle} satisfies \eqref{eq:conv0}. Hence, the result of Theorem \ref{main} applies for $\tilde\theta_n(\bold{X}_n) = \hat\theta_n(\bold{X}_n)$.
\end{remark}

\begin{remark}
In Theorem \ref{main}  there is no need to impose $C_\phi=0$ and $K_\phi=1$ as, e.g. in Morales {\it et al.} (1997). Of course, in our case the constants $C_\phi$ and $K_\phi$ enter in the asymptotic distribution of the test statistics.
The convergence result is also interesting because, contrary to the i.i.d case,  the rate of convergence of the estimators of $\theta$ for the drift and diffusion coefficients are different and are respectively equal to $\sqrt{n\Delta_n}$ and $\sqrt{n}$.
\end{remark}

\begin{remark}
As remarked in Uchida and Yoshida (2001), it is always better to derive approximate ML estimators and the test statistics on different approximations of the true likelihood to avoid circularities.
\end{remark}

\section{Numerical analysis}\label{sec:applications}
Although asymptotic properties have been obtained, what really
matters in application is the behaviour of the test statistics
under fine sample setup. We study the empirical performance of the
test for small samples in terms of level of the test and power
under some alternatives. In the analysis we consider the estimator
\eqref{eq:mle} and the following quantities
\begin{itemize}
\item estimated $\alpha$-divergences
$$
\mathbb{D}_\alpha(\hat\theta_n(\bold{X}_n),\theta_0) =
\phi_\alpha\left(\frac{f_n({\bf X}_n,\hat\theta_n({\bf
X}_n))}{f_n({\bf X}_n,\theta_0)}\right)
$$
with $\phi_\alpha(x) = 4(1-x^\frac{1+\alpha}{2})/(1-\alpha^2)$,
with $C_{\alpha} = \frac{2}{\alpha-1}$ and $K_\phi=1$. We consider
$\alpha\in\{-0.99, -0.90, -0.75, -0.50, -0.25, -0.10\}$; \item
estimated power-divergences
$$
\mathbb{D}_\lambda(\hat\theta_n(\bold{X}_n),\theta_0) =
\phi_\lambda\left(\frac{f_n({\bf X}_n,\hat\theta_n({\bf
X}_n))}{f_n({\bf X}_n,\theta_0)}\right)
$$
with $\phi_\lambda(x) = (x^{\lambda+1} - x -
\lambda(x-1))/(\lambda(\lambda+1))$, with $C_\lambda=0$,
$K_\lambda=1$. We consider $\lambda \in\{-0.99, -1.20, -1.50,
-1.75, -2.00, -2.50\}$;
\item likelihood ratio statistic
$$\mathbb{D}_{\log}(\hat\theta_n(\bold{X}_n),\theta_0)  = -\log\left(\frac{f_n({\bf X}_n,\hat\theta_n({\bf X}_n))}{f_n({\bf X}_n,\theta_0)}\right)
$$
\end{itemize}
For $\mathbb D_\alpha$ and $\mathbb D_\lambda$, the threshold of
the rejection region of the test  are calculated using formula
\eqref{chi3} as the empirical quantiles of  \eqref{chi3} of 100000
simulations of the random variable $\chi^2_{p+q}$. For
$\mathbb{D}_{\log}$ is again used formula  \eqref{chi3} but exact
quantiles of the random variable $\chi^2_{p+q}$ are used. Because
the interest is in testing $\mathbb D_\phi = 0$ against $\mathbb
D_\phi \neq 0$, whenever $f_n({\bf X}_n,\tilde\theta_n({\bf
X}_n))> f_n({\bf X}_n,\theta_0)$ we exchange the numerator and the
denominator to avoid negative signs in the test statistics.
Usually, this is not going to happen if $\phi$ is convex and
$\phi'(1)=0$ (see, e.g. Morales {\it et al.}, 1997).

We evaluate the empirical level of the test calculated as the number of times the test rejects the null hypothesis under the true model, i.e.
$$
\hat\alpha_n = \frac{1}{M} \sum_{i=1}^M {\bf 1}_{\{\mathbb D_\phi
> c_{\alpha}\}}
$$
where ${\bf 1}_A$ is the indicator function of set $A$, $M=10000$
is the number of simulations and $c_{\alpha}$ is the
$(1-\alpha)\%$ quantile of the proper distribution. Similarly we
calculate the power of the test under alternative models as
$$
\hat\beta_n = \frac{1}{M} \sum_{i=1}^M {\bf 1}_{\{\mathbb D_\phi >
c_{\alpha}\}}
$$
In our experiments we consider the two families of stochastic processes borrowed from finance
\begin{itemize}
\item the Vasicek (VAS) model
$$
\de X_t = \kappa(\alpha-X_t)\de t + \sigma X_t \de W_t
$$
where, in finance, $\sigma$ is interpreted as volatility, $\alpha$ is the long-run equilibrium value of the process and $\kappa$ is the speed of reversion. Let $(\kappa_0, \alpha_0, \sigma^2_0) = (0.85837, 0.089102, 0.0021854)$,
we consider three different sets of hypotheses for the parameters
\begin{center}
\begin{tabular}{lc}
model & $\theta = (\kappa, \alpha, \sigma^2)$\\
\hline
$\text{VAS}_0$ &   $(\kappa_0, \alpha_0, \sigma^2_0)$\\
$\text{VAS}_1$ &   $(4 \cdot \kappa_0, \alpha_0, 4 \cdot \sigma^2_0 )$\\
$\text{VAS}_2$ &   $(\frac14 \kappa_0, \alpha_0, \frac14 \cdot \sigma^2_0)$
\end{tabular}
\end{center}
The interesting facts are that $\text{VAS}_0$, $\text{VAS}_1$ and $\text{VAS}_2$ have all the same
stationary distributions $N(\alpha_0, \sigma_0^2/(2\kappa_0))$, a Gaussian transition density
$$N\left(\alpha_0+(x_0-\alpha_0)e^{-\kappa t},  \frac{\sigma_0^2 (1-e^{-2\kappa t})}{2\kappa_0} \right)
$$
and covariance function given by
$$
{\rm Cov}(X_s, X_t) = \frac{\sigma_0^2}{2\kappa_0} e^{-\kappa(s+t)}\left(
e^{-2\kappa(s\wedge t)-1}\right)
$$
and both show a strong dependency of the covariance as a function
of $\kappa$, which makes this model interesting in comparison with
the i.i.d. setting; \item the Cox-Ingersoll-Ross (CIR) model
$$
\de X_t = \kappa(\alpha-X_t)\de t + \sigma \sqrt{X_t} \de W_t
$$
Let $(\kappa_0, \alpha_0, \sigma^2_0) = (0.89218, 0.09045, 0.032742)$,
we consider different sets of hypotheses for the parameters
\begin{center}
\begin{tabular}{lc}
model & $\theta = (\kappa, \alpha, \sigma^2)$\\
\hline
$\text{CIR}_0$ &   $(\kappa_0, \alpha_0, \sigma_0^2)$\\
$\text{CIR}_1$ &   $(\frac12 \cdot \kappa_0, \alpha_0,\frac12 \cdot  \sigma_0^2)$\\
$\text{CIR}_2$ &   $(\frac14 \cdot \kappa_0, \alpha_0,\frac14 \cdot  \sigma_0^2)$
\end{tabular}
\end{center}
This model has a transition density of $\chi^2$-type, hence local gaussian approximation is less likely to hold for non negligible values of $\Delta_n$.
\end{itemize}
The parameters of the above models, have been chosen
according to Pritsker (1998) and Chen {\it et al.} (2008), in particular  $\text{VAS}_0$ corresponds to the model estimated by A\"it-Sahalia (1996) for real interest rates data.

We study the level and the power of the three family of test statistics for different values
of $\Delta_n \in \{0.1, 0.001\}$ and $n \in \{50, 100, 500\}$.
For the same trajectory, hence we simulate 1000 observations and we extract
only that last $n$ observations. Disregarding the first part of the trajectory ensures that the process
is in the stationary state.

The results of these simulations are reported in the Tables
\ref{tabvas13l}-\ref{tabcir3a}. We point out that in the Tables
\ref{tabvas1a}, \ref{tabvas3a}, \ref{tabcir1a} and \ref{tabcir3a},
in the column ``model ($\alpha$, $n$)'' the $\alpha$ corresponds
to the true level of the test used to calculate $c_{\alpha}$. The
other $\alpha$'s in the first row of the tables correspond to the
$\alpha$ in $\phi_\alpha$-divergences.

\paragraph{Summary of the analysis for the Vasicek model}
It turns out that $\alpha$-divergences are not very good in terms of estimated
level of the test, but their power function behaves as expected. It also emerges that
 for $\lambda=-0.99$, the power divergence cannot identify as wrong model VAS$_1$ for small
 sample size $n=50$ and $\Delta_n=0.001$ (Table \ref{tabvas1p}, row 2), although this is not the case for the power-divergences and the likelihood ratio test (Tables \ref{tabvas1a} and \ref{tabvas13l}, row 2).

In general power divergences for $\lambda$ in $\{-0.99, -1.20,
-1.50, -1.75, -2.00\}$ have always very small estimated level and
high power under the selected alternatives. The
$\alpha$-divergences, do not behave very good and, the way they
are defined, only approximate the likelihood ratio for
$\alpha=-0.99$.

The power divergences are, on average, better than the likelihood ratio test in terms of both empirical level $\hat\alpha$ and power $\hat\beta$ under the selected alternatives.

\paragraph{Summary of the analysis for the CIR model}
The same average considerations apply to the case of CIR model.
The difference is that, for small sample size, all test statistics
have low power under the alternative CIR$_1$ while CIR$_2$ doesn't
present particular problems.

\section{Proofs}\label{sec:proofs}
The following important Lemmas are useful to prove the Theorem
\ref{main}.

\begin{lemma}[Kessler, 1997]\label{lem3} Under the assumptions \ref{ass:1}-\ref{ass:3}, as $n\Delta_n^2\to 0$ the following hold true

\begin{equation}
\label{eq2} \Gamma^\frac12 \nabla_\theta g_n({\bf X}_n, \theta_0)
\stackrel{p}{\to} N(0, \mathcal I(\theta_0))
\end{equation}
\end{lemma}
\begin{lemma}[Uchida and Yoshida, 2005]\label{lem4} Under the assumptions \ref{ass:1}-\ref{ass:3}, as $n\Delta_n^2\to 0$ the following hold  true
\begin{equation}
\label{eq3} \Gamma^\frac12\nabla_\theta l_n({\bf X}_n, \theta_0)
=\Gamma^\frac12\nabla_\theta g_n({\bf X}_n, \theta_0)+o_p(1)
\end{equation}
\end{lemma}
\begin{lemma}[Uchida and Yoshida, 2005]\label{lem5} Under the assumptions \ref{ass:1}-\ref{ass:3}, as $n\Delta_n^2\to 0$ the following hold  true
\begin{equation}
\label{eq31} \Gamma^\frac12\nabla_\theta^2 l_n({\bf X}_n, \theta_0)
\Gamma^\frac12\stackrel{p}{\to}-\mathcal{I}(\theta_0)
\end{equation}
\end{lemma}

\begin{proof}[Proof of Theorem \ref{main}]
We start by applying {\it delta} method. We denote the gradient
vector by $\nabla_\theta = [\partial/\partial\theta_i]$, $i=1,
\ldots, p+q$ and similarly the Hessian matrix by  $\nabla^2_\theta
= [\partial^2/\partial\theta_i\partial\theta_j]$, $i,j =1, \ldots,
p+q$.

$i)$ We can write that
\begin{eqnarray*}
\mathbb{D}_\phi(\tilde\theta_n(\bold{X}_n),\theta_0)&=&\mathbb{D}_\phi(\theta_0,\theta_0)+[\nabla_\theta
\mathbb{D}_\phi(\theta_0,\theta_0)]^T(\tilde\theta_n(\bold{X}_n)-\theta_0)+o_p(1)\\
&=&[\nabla_\theta
\mathbb{D}_\phi(\theta_0,\theta_0)]^T(\tilde\theta_n(\bold{X}_n)-\theta_0)+o_p(1)
\end{eqnarray*}
because $\mathbb{D}_\phi(\theta_0,\theta_0)=0$. Noting that for
$k=1,...,p+q$
$$
\frac{\partial}{\partial\theta_k}\left[\phi\left(\frac{f_n(\cdot,\theta)}{f_n(\cdot,\theta_0)}\right)
\right]=\frac{1}{f_n(\cdot,\theta_0)}\phi'\left(\frac{f_n(\cdot,\theta)}{f_n(\cdot,\theta_0)}\right)
\frac{\partial f_n(\cdot,\theta)}{\partial\theta_k}
$$
 by Assumption \ref{ass:C1} follows that
\begin{eqnarray*}
\nabla_\theta \mathbb{D}_\phi(\theta_0,\theta_0)=C_\phi\left.
\nabla_\theta
l_n(\bold{X}_n,\theta)\right|_{\theta=\theta_0}=C_\phi\nabla_\theta
l_n(\bold{X}_n,\theta_0)\end{eqnarray*}
and therefore
\begin{eqnarray}\label{eq:proof1}
\mathbb{D}_\phi(\tilde\theta_n(\bold{X}_n),\theta_0)=
C_\phi\left[\Gamma^\frac12\nabla_\theta
l_n(\bold{X}_n,\theta_0)\right]^T\Gamma^{-\frac12}(\tilde\theta_n(\bold{X}_n)-\theta_0)+o_p(1)
\end{eqnarray}
From \eqref{eq:proof1} by means of Lemma \ref{lem4}-\ref{lem3} and
Slutsky's Theorem immediately follows
$$\mathbb{D}_\phi(\tilde\theta_n(\bold{X}_n),\theta_0)\stackrel{d}{\to}C_\phi\chi_{p+q}^2$$

$ii)$ Since for $k,j=1,...,p+q$
\begin{eqnarray*}
\frac{\partial^2}{\partial\theta_k\partial\theta_j}\left[\phi\left(\frac{f_n(\cdot,\theta)}{f_n(\cdot,\theta_0)}\right)
\right]&=&\frac{1}{f_n^2(\cdot,\theta_0)}\phi''\left(\frac{f_n(\cdot,\theta)}{f_n(\cdot,\theta_0)}\right)
\frac{\partial f_n(\cdot,\theta)}{\partial\theta_k}
\frac{\partial f_n(\cdot,\theta)}{\partial\theta_j}\\
&&+\frac{1}{f_n(\cdot,\theta_0)}\phi'\left(\frac{f_n(\cdot,\theta)}{f_n(\cdot,\theta_0)}\right)
\frac{\partial^2
f_n(\cdot,\theta)}{\partial\theta_k\partial\theta_j}
\end{eqnarray*}
follows that
\begin{eqnarray*}
\mathbb{D}_\phi(\tilde\theta_n(\bold{X}_n),\theta_0)&=&\frac{1}{2}[\Gamma^{-1/2}(\tilde\theta_n(\bold{X}_n)-\theta_0)]^T\Gamma^{1/2}\nabla_\theta^2\mathbb{D}_\phi(\theta_0,\theta_0)\Gamma^{1/2}\\
&&\times \Gamma^{-1/2}(\tilde\theta_n(\bold{X}_n)-\theta_0)+o_p(1)\\
&=&\frac{K_\phi}{2}[\Gamma^{-1/2}(\tilde\theta_n(\bold{X}_n)-\theta_0)]^T\Gamma^{1/2}\nabla_\theta
l_n(\bold{X}_n,\theta_0)[\Gamma^{1/2}\nabla_\theta
l_n(\bold{X}_n,\theta_0)]^T\\
&&\times \Gamma^{-1/2}(\tilde\theta_n(\bold{X}_n)-\theta_0)+o_p(1)
\end{eqnarray*}
From \eqref{eq:proof1} by means of Lemma \ref{lem4}-\ref{lem3} and
Slutsky's Theorem immediately follows

$$\mathbb{D}_\phi(\tilde\theta_n(\bold{X}_n),\theta_0)\stackrel{d}{\to}
\frac{K_\phi}{2}Z_{p+q}$$ It's easy to verify that the density
function of the r.v. $Z_{p+q}$ is equal to
\begin{equation}
f_{Z_{p+q}}(z)=\frac{(1/2)^{\frac{p+q}{2}}}{\Gamma\left(\frac{p+q}{2}\right)}
\sqrt{z}^{\frac{p+q}{2}-1}e^{-\sqrt{z}/2}\frac{1}{2\sqrt{z}},\quad
z> 0
\end{equation}

 $iii)$ By previous considerations we have that
 \begin{eqnarray}\label{eq:proof3}
\mathbb{D}_\phi(\tilde\theta_n(\bold{X}_n),\theta_0)&=&[\nabla_\theta
\mathbb{D}_\phi(\theta_0,\theta_0)]^T(\tilde\theta_n(\bold{X}_n)-\theta_0)\notag\\
&&+\frac{1}{2}[\Gamma^{-1/2}(\tilde\theta_n(\bold{X}_n)-\theta_0)]^T\Gamma^{1/2}\nabla_\theta^2\mathbb{D}_\phi(\theta_0,\theta_0)\Gamma^{1/2}\notag\\
&&\times \Gamma^{-1/2}(\tilde\theta_n(\bold{X}_n)-\theta_0)+o_p(1)
\end{eqnarray}
where
 \begin{eqnarray}\label{eq:proof4}
&&\nabla_\theta^2\mathbb{D}_\phi(\theta_0,\theta_0)\notag\\
&&=K_\phi\nabla_\theta l_n(\bold{X}_n,\theta_0)[\nabla_\theta
l_n(\bold{X}_n,\theta_0)]^T+C_\phi\frac{1}{f(\bold{X}_n,\theta_0)}\nabla_\theta^2f(\bold{X}_n,\theta_0)
\notag\\
&&=(K_\phi+C_\phi) \nabla_\theta
l_n(\bold{X}_n,\theta_0)[\nabla_\theta
l_n(\bold{X}_n,\theta_0)]^T+C_\phi\nabla_\theta^2l_n({\bf X}_n,
\theta_0)
\end{eqnarray}
Plugging in \eqref{eq:proof3} the quantity \eqref{eq:proof4} we
derive, applying again Lemma \ref{lem3}-\ref{lem5}, the following
result
\begin{equation*}
\mathbb{D}_\phi(\tilde\theta_n(\bold{X}_n),\theta_0)\stackrel{d}{\to}
\frac12\left[C_\phi\chi^2_{p+q}+(C_\phi+K_\phi)Z_{p+q}\right]
\end{equation*}

\end{proof}

\section*{Conclusions}
It seems that, as in the i.i.d. case, also for discretely observed diffusion processes the $\phi$-divergences may compete or improve the performance of the standard  likelihood ratio statistics.
In particular, the power divergences are in general quite good in terms of estimated level and power of the test even for moderate sample sizes (e.g. $n\geq 100$ in our simulations).

The package {\tt sde} for the {\tt R} statistical environment (R Development Core Team, 2008) and freely available at \url{http://cran.R-Project.org} contains the function {\tt sdeDiv} which implements the  $\phi$-divergence test statistics.

\newpage

\begin{table}
\begin{center}
\begin{tabular}{lrr}
  \hline
model (n) & $\alpha=0.01$ & $\alpha=0.05$ \\
  \hline
VAS$_0$ (50) & 0.01 & 0.04 \\
  VAS$_1$ (50) & 1.00 & 1.00 \\
  VAS$_2$ (50) & 1.00 & 1.00 \\
   &  &  \\
  VAS$_0$ (100) & 0.01 & 0.04 \\
  VAS$_1$ (100) & 1.00 & 1.00 \\
  VAS$_2$ (100) & 1.00 & 1.00 \\
   &  &  \\
  VAS$_0$ (500) & 0.01 & 0.07 \\
  VAS$_1$ (500) & 1.00 & 1.00 \\
  VAS$_2$ (500) & 1.00 & 1.00 \\
   \hline
\end{tabular}\\
\begin{tabular}{lrr}
  \hline
model (n) & $\alpha=0.01$ & $\alpha=0.05$ \\
  \hline
VAS$_0$ (50) & 0.01 & 0.04 \\
  VAS$_1$ (50) & 1.00 & 1.00 \\
  VAS$_2$ (50) & 1.00 & 1.00 \\
   &  &  \\
  VAS$_0$ (100) & 0.01 & 0.04 \\
  VAS$_1$ (100) & 1.00 & 1.00 \\
  VAS$_2$ (100) & 1.00 & 1.00 \\
   &  &  \\
  VAS$_0$ (500) & 0.00 & 0.02 \\
  VAS$_1$ (500) & 1.00 & 1.00 \\
  VAS$_2$ (500) & 1.00 & 1.00 \\
   \hline
\end{tabular}
\end{center}
\caption{Numbers represent probability of rejection under the true generating
model, with $c_\alpha$ calculated under $H_0$. Therefore, the values are
 $\hat\alpha$ under model ``$0$'' and $\hat\beta$ otherwise. Estimates
 calculated on 10000 experiments.
Likelihood
ratio, for $\Delta_n=0.001$ (up) and $\Delta_n=0.1$ (bottom).}
\label{tabvas13l}
\end{table}

\begin{landscape}
\begin{table}
\begin{center}
\begin{tabular}{lrrrrrr}
  \hline
model ($\alpha$, $n$) & $\alpha=-0.99$ & $\alpha=-0.90$ & $\alpha=-0.75$ & $\alpha=-0.50$ & $\alpha=-0.25$ & $\alpha=-0.10$ \\
  \hline
VAS$_0$ (0.01, 50) & 0.01 & 0.10 & 0.39 & 0.62 & 0.73 & 0.77 \\
  VAS$_1$ (0.01, 50) & 1.00 & 1.00 & 1.00 & 1.00 & 1.00 & 1.00 \\
  VAS$_2$ (0.01, 50) & 1.00 & 1.00 & 1.00 & 1.00 & 1.00 & 1.00 \\
   &  &  &  &  &  &  \\
  VAS$_0$ (0.05, 50) & 0.04 & 0.12 & 0.39 & 0.62 & 0.73 & 0.77 \\
  VAS$_1$ (0.05, 50) & 1.00 & 1.00 & 1.00 & 1.00 & 1.00 & 1.00 \\
  VAS$_2$ (0.05, 50) & 1.00 & 1.00 & 1.00 & 1.00 & 1.00 & 1.00 \\
   &  &  &  &  &  &  \\
   &  &  &  &  &  &  \\
  VAS$_0$ (0.01, 100) & 0.01 & 0.10 & 0.39 & 0.63 & 0.74 & 0.78 \\
  VAS$_1$ (0.01, 100) & 1.00 & 1.00 & 1.00 & 1.00 & 1.00 & 1.00 \\
  VAS$_2$ (0.01, 100) & 1.00 & 1.00 & 1.00 & 1.00 & 1.00 & 1.00 \\
   &  &  &  &  &  &  \\
  VAS$_0$ (0.05, 100) & 0.04 & 0.11 & 0.40 & 0.63 & 0.74 & 0.78 \\
  VAS$_1$ (0.05, 100) & 1.00 & 1.00 & 1.00 & 1.00 & 1.00 & 1.00 \\
  VAS$_2$ (0.05, 100) & 1.00 & 1.00 & 1.00 & 1.00 & 1.00 & 1.00 \\
   &  &  &  &  &  &  \\
   &  &  &  &  &  &  \\
  VAS$_0$ (0.01, 500) & 0.02 & 0.18 & 0.61 & 0.83 & 0.90 & 0.92 \\
  VAS$_1$ (0.01, 500) & 1.00 & 1.00 & 1.00 & 1.00 & 1.00 & 1.00 \\
  VAS$_2$ (0.01, 500) & 1.00 & 1.00 & 1.00 & 1.00 & 1.00 & 1.00 \\
   &  &  &  &  &  &  \\
  VAS$_0$ (0.05, 500) & 0.07 & 0.20 & 0.61 & 0.83 & 0.90 & 0.92 \\
  VAS$_1$ (0.05, 500) & 1.00 & 1.00 & 1.00 & 1.00 & 1.00 & 1.00 \\
  VAS$_2$ (0.05, 500) & 1.00 & 1.00 & 1.00 & 1.00 & 1.00 & 1.00 \\
   \hline
\end{tabular}
\end{center}
\caption{Numbers represent probability of rejection under the true generating
model, with $c_\alpha$ calculated under $H_0$. Therefore, the values are
 $\hat\alpha$ under model ``$0$'' and $\hat\beta$ otherwise. Estimates
 calculated on 10000 experiments.
 $\alpha$-divergences, for $\Delta_n=0.001$.}
\label{tabvas1a}
\end{table}

\begin{table}
\begin{center}
\begin{tabular}{lrrrrrr}
  \hline
model ($\alpha$, $n$) & $\lambda=-0.99$ & $\lambda=-1.20$ & $\lambda=-1.50$ & $\lambda=-1.75$ & $\lambda=-2.00$ & $\lambda=-2.50$ \\
  \hline
VAS$_0$ (0.01, 50) & 0.00 & 0.00 & 0.00 & 0.01 & 0.02 & 0.04 \\
  VAS$_1$ (0.01, 50) & 0.00 & 0.99 & 1.00 & 1.00 & 1.00 & 1.00 \\
  VAS$_2$ (0.01, 50) & 0.40 & 1.00 & 1.00 & 1.00 & 1.00 & 1.00 \\
   &  &  &  &  &  &  \\
  VAS$_0$ (0.05, 50) & 0.00 & 0.00 & 0.00 & 0.01 & 0.03 & 0.06 \\
  VAS$_1$ (0.05, 50) & 0.67 & 1.00 & 1.00 & 1.00 & 1.00 & 1.00 \\
  VAS$_2$ (0.05, 50) & 0.99 & 1.00 & 1.00 & 1.00 & 1.00 & 1.00 \\
   &  &  &  &  &  &  \\
   &  &  &  &  &  &  \\
  VAS$_0$ (0.01, 100) & 0.00 & 0.00 & 0.00 & 0.01 & 0.02 & 0.04 \\
  VAS$_1$ (0.01, 100) & 0.23 & 1.00 & 1.00 & 1.00 & 1.00 & 1.00 \\
  VAS$_2$ (0.01, 100) & 0.88 & 1.00 & 1.00 & 1.00 & 1.00 & 1.00 \\
   &  &  &  &  &  &  \\
  VAS$_0$ (0.05, 100) & 0.00 & 0.00 & 0.00 & 0.01 & 0.03 & 0.06 \\
  VAS$_1$ (0.05, 100) & 1.00 & 1.00 & 1.00 & 1.00 & 1.00 & 1.00 \\
  VAS$_2$ (0.05, 100) & 1.00 & 1.00 & 1.00 & 1.00 & 1.00 & 1.00 \\
   &  &  &  &  &  &  \\
   &  &  &  &  &  &  \\
  VAS$_0$ (0.01, 500) & 0.00 & 0.00 & 0.00 & 0.01 & 0.03 & 0.08 \\
  VAS$_1$ (0.01, 500) & 1.00 & 1.00 & 1.00 & 1.00 & 1.00 & 1.00 \\
  VAS$_2$ (0.01, 500) & 1.00 & 1.00 & 1.00 & 1.00 & 1.00 & 1.00 \\
   &  &  &  &  &  &  \\
  VAS$_0$ (0.05, 500) & 0.00 & 0.00 & 0.01 & 0.03 & 0.06 & 0.12 \\
  VAS$_1$ (0.05, 500) & 1.00 & 1.00 & 1.00 & 1.00 & 1.00 & 1.00 \\
  VAS$_2$ (0.05, 500) & 1.00 & 1.00 & 1.00 & 1.00 & 1.00 & 1.00 \\
   \hline
\end{tabular}
\end{center}
\caption{Numbers represent probability of rejection under the true generating
model, with $c_\alpha$ calculated under $H_0$. Therefore, the values are
 $\hat\alpha$ under model ``$0$'' and $\hat\beta$ otherwise. Estimates
 calculated on 10000 experiments.
Power-divergences  for $\Delta_n=0.001$} \label{tabvas1p}
\end{table}

\begin{table}
\begin{center}
\begin{tabular}{lrrrrrr}
  \hline
model ($\alpha$, $n$) & $\alpha=-0.99$ & $\alpha=-0.90$ & $\alpha=-0.75$ & $\alpha=-0.50$ & $\alpha=-0.25$ & $\alpha=-0.10$ \\
  \hline
VAS$_0$ (0.01, 50) & 0.01 & 0.15 & 0.55 & 0.78 & 0.86 & 0.88 \\
  VAS$_1$ (0.01, 50) & 1.00 & 1.00 & 1.00 & 1.00 & 1.00 & 1.00 \\
  VAS$_2$ (0.01, 50) & 1.00 & 1.00 & 1.00 & 1.00 & 1.00 & 1.00 \\
   &  &  &  &  &  &  \\
  VAS$_0$ (0.05, 50) & 0.05 & 0.17 & 0.55 & 0.78 & 0.86 & 0.88 \\
  VAS$_1$ (0.05, 50) & 1.00 & 1.00 & 1.00 & 1.00 & 1.00 & 1.00 \\
  VAS$_2$ (0.05, 50) & 1.00 & 1.00 & 1.00 & 1.00 & 1.00 & 1.00 \\
   &  &  &  &  &  &  \\
   &  &  &  &  &  &  \\
  VAS$_0$ (0.01, 100) & 0.01 & 0.13 & 0.48 & 0.71 & 0.80 & 0.83 \\
  VAS$_1$ (0.01, 100) & 1.00 & 1.00 & 1.00 & 1.00 & 1.00 & 1.00 \\
  VAS$_2$ (0.01, 100) & 1.00 & 1.00 & 1.00 & 1.00 & 1.00 & 1.00 \\
   &  &  &  &  &  &  \\
  VAS$_0$ (0.05, 100) & 0.04 & 0.15 & 0.48 & 0.71 & 0.80 & 0.83 \\
  VAS$_1$ (0.05, 100) & 1.00 & 1.00 & 1.00 & 1.00 & 1.00 & 1.00 \\
  VAS$_2$ (0.05, 100) & 1.00 & 1.00 & 1.00 & 1.00 & 1.00 & 1.00 \\
   &  &  &  &  &  &  \\
   &  &  &  &  &  &  \\
  VAS$_0$ (0.01, 500) & 0.00 & 0.06 & 0.25 & 0.54 & 0.69 & 0.74 \\
  VAS$_1$ (0.01, 500) & 1.00 & 1.00 & 1.00 & 1.00 & 1.00 & 1.00 \\
  VAS$_2$ (0.01, 500) & 1.00 & 1.00 & 1.00 & 1.00 & 1.00 & 1.00 \\
   &  &  &  &  &  &  \\
  VAS$_0$ (0.05, 500) & 0.02 & 0.07 & 0.25 & 0.54 & 0.69 & 0.74 \\
  VAS$_1$ (0.05, 500) & 1.00 & 1.00 & 1.00 & 1.00 & 1.00 & 1.00 \\
  VAS$_2$ (0.05, 500) & 1.00 & 1.00 & 1.00 & 1.00 & 1.00 & 1.00 \\
   \hline
\end{tabular}
\end{center}
\caption{Numbers represent probability of rejection under the true generating
model, with $c_\alpha$ calculated under $H_0$. Therefore, the values are
 $\hat\alpha$ under model ``$0$'' and $\hat\beta$ otherwise. Estimates
 calculated on 10000 experiments. $\alpha$-divergences, for $\Delta_n=0.1$}
\label{tabvas3a}
\end{table}

\begin{table}
\begin{center}
\begin{tabular}{lrrrrrr}
  \hline
model ($\alpha$, $n$) & $\lambda=-0.99$ & $\lambda=-1.20$ & $\lambda=-1.50$ & $\lambda=-1.75$ & $\lambda=-2.00$ & $\lambda=-2.50$ \\
  \hline
VAS$_0$ (0.01, 50) & 0.00 & 0.00 & 0.00 & 0.01 & 0.02 & 0.05 \\
  VAS$_1$ (0.01, 50) & 1.00 & 1.00 & 1.00 & 1.00 & 1.00 & 1.00 \\
  VAS$_2$ (0.01, 50) & 1.00 & 1.00 & 1.00 & 1.00 & 1.00 & 1.00 \\
   &  &  &  &  &  &  \\
  VAS$_0$ (0.05, 50) & 0.00 & 0.00 & 0.00 & 0.02 & 0.03 & 0.09 \\
  VAS$_1$ (0.05, 50) & 1.00 & 1.00 & 1.00 & 1.00 & 1.00 & 1.00 \\
  VAS$_2$ (0.05, 50) & 1.00 & 1.00 & 1.00 & 1.00 & 1.00 & 1.00 \\
   &  &  &  &  &  &  \\
   &  &  &  &  &  &  \\
  VAS$_0$ (0.01, 100) & 0.00 & 0.00 & 0.00 & 0.00 & 0.01 & 0.05 \\
  VAS$_1$ (0.01, 100) & 1.00 & 1.00 & 1.00 & 1.00 & 1.00 & 1.00 \\
  VAS$_2$ (0.01, 100) & 1.00 & 1.00 & 1.00 & 1.00 & 1.00 & 1.00 \\
   &  &  &  &  &  &  \\
  VAS$_0$ (0.05, 100) & 0.00 & 0.00 & 0.00 & 0.01 & 0.03 & 0.08 \\
  VAS$_1$ (0.05, 100) & 1.00 & 1.00 & 1.00 & 1.00 & 1.00 & 1.00 \\
  VAS$_2$ (0.05, 100) & 1.00 & 1.00 & 1.00 & 1.00 & 1.00 & 1.00 \\
   &  &  &  &  &  &  \\
   &  &  &  &  &  &  \\
  VAS$_0$ (0.01, 500) & 0.00 & 0.00 & 0.00 & 0.00 & 0.01 & 0.02 \\
  VAS$_1$ (0.01, 500) & 1.00 & 1.00 & 1.00 & 1.00 & 1.00 & 1.00 \\
  VAS$_2$ (0.01, 500) & 1.00 & 1.00 & 1.00 & 1.00 & 1.00 & 1.00 \\
   &  &  &  &  &  &  \\
  VAS$_0$ (0.05, 500) & 0.00 & 0.00 & 0.00 & 0.01 & 0.01 & 0.04 \\
  VAS$_1$ (0.05, 500) & 1.00 & 1.00 & 1.00 & 1.00 & 1.00 & 1.00 \\
  VAS$_2$ (0.05, 500) & 1.00 & 1.00 & 1.00 & 1.00 & 1.00 & 1.00 \\
   \hline
\end{tabular}
\end{center}
\caption{Numbers represent probability of rejection under the true generating
model, with $c_\alpha$ calculated under $H_0$. Therefore, the values are
 $\hat\alpha$ under model ``$0$'' and $\hat\beta$ otherwise. Estimates
 calculated on 10000 experiments.
Power-divergences  for $\Delta_n=0.1$} \label{tabvas3p}
\end{table}
\end{landscape}

\begin{table}
\begin{center}
\begin{tabular}{lrr}
  \hline
model (n) & $\alpha=0.01$ & $\alpha=0.05$ \\
  \hline
CIR$_0$ (50) & 0.02 & 0.11 \\
  CIR$_1$ (50) & 0.59 & 0.84 \\
  CIR$_2$ (50) & 1.00 & 1.00 \\
   &  &  \\
  CIR$_0$ (100) & 0.03 & 0.11 \\
  CIR$_1$ (100) & 0.96 & 0.99 \\
  CIR$_2$ (100) & 1.00 & 1.00 \\
   &  &  \\
  CIR$_0$ (500) & 0.02 & 0.09 \\
  CIR$_1$ (500) & 1.00 & 1.00 \\
  CIR$_2$ (500) & 1.00 & 1.00 \\
   \hline
\end{tabular}\\
\begin{tabular}{lrr}
  \hline
model (n) & $\alpha=0.01$ & $\alpha=0.05$ \\
  \hline
CIR$_0$ (50) & 0.01 & 0.04 \\
  CIR$_1$ (50) & 0.78 & 0.93 \\
  CIR$_2$ (50) & 1.00 & 1.00 \\
   &  &  \\
  CIR$_0$ (100) & 0.01 & 0.04 \\
  CIR$_1$ (100) & 0.99 & 1.00 \\
  CIR$_2$ (100) & 1.00 & 1.00 \\
   &  &  \\
  CIR$_0$ (500) & 0.00 & 0.02 \\
  CIR$_1$ (500) & 1.00 & 1.00 \\
  CIR$_2$ (500) & 1.00 & 1.00 \\
   \hline
\end{tabular}
\end{center}
\caption{Numbers represent probability of rejection under the true
model, with rejection region calculated under $H_0$. Likelihood
ratio, for $\Delta_n=0.001$ (up) and $\Delta_n=0.1$ (bottom).}
\label{tabcir13l}
\end{table}

\begin{landscape}
\begin{table}
\begin{center}
\begin{tabular}{lrrrrrr}
  \hline
model ($\alpha$, $n$) & $\alpha=-0.99$ & $\alpha=-0.90$ & $\alpha=-0.75$ & $\alpha=-0.50$ & $\alpha=-0.25$ & $\alpha=-0.10$ \\
  \hline
CIR$_0$ (0.01, 50) & 0.03 & 0.28 & 0.69 & 0.83 & 0.89 & 0.90 \\
  CIR$_1$ (0.01, 50) & 0.63 & 0.95 & 1.00 & 1.00 & 1.00 & 1.00 \\
  CIR$_2$ (0.01, 50) & 1.00 & 1.00 & 1.00 & 1.00 & 1.00 & 1.00 \\
   &  &  &  &  &  &  \\
  CIR$_0$ (0.05, 50) & 0.12 & 0.31 & 0.69 & 0.83 & 0.89 & 0.90 \\
  CIR$_1$ (0.05, 50) & 0.86 & 0.96 & 1.00 & 1.00 & 1.00 & 1.00 \\
  CIR$_2$ (0.05, 50) & 1.00 & 1.00 & 1.00 & 1.00 & 1.00 & 1.00 \\
   &  &  &  &  &  &  \\
   &  &  &  &  &  &  \\
  CIR$_0$ (0.01, 100) & 0.03 & 0.28 & 0.69 & 0.85 & 0.89 & 0.91 \\
  CIR$_1$ (0.01, 100) & 0.97 & 1.00 & 1.00 & 1.00 & 1.00 & 1.00 \\
  CIR$_2$ (0.01, 100) & 1.00 & 1.00 & 1.00 & 1.00 & 1.00 & 1.00 \\
   &  &  &  &  &  &  \\
  CIR$_0$ (0.05, 100) & 0.12 & 0.31 & 0.69 & 0.85 & 0.89 & 0.91 \\
  CIR$_1$ (0.05, 100) & 0.99 & 1.00 & 1.00 & 1.00 & 1.00 & 1.00 \\
  CIR$_2$ (0.05, 100) & 1.00 & 1.00 & 1.00 & 1.00 & 1.00 & 1.00 \\
   &  &  &  &  &  &  \\
   &  &  &  &  &  &  \\
  CIR$_0$ (0.01, 500) & 0.03 & 0.22 & 0.59 & 0.79 & 0.86 & 0.88 \\
  CIR$_1$ (0.01, 500) & 1.00 & 1.00 & 1.00 & 1.00 & 1.00 & 1.00 \\
  CIR$_2$ (0.01, 500) & 1.00 & 1.00 & 1.00 & 1.00 & 1.00 & 1.00 \\
   &  &  &  &  &  &  \\
  CIR$_0$ (0.05, 500) & 0.09 & 0.24 & 0.59 & 0.79 & 0.86 & 0.89 \\
  CIR$_1$ (0.05, 500) & 1.00 & 1.00 & 1.00 & 1.00 & 1.00 & 1.00 \\
  CIR$_2$ (0.05, 500) & 1.00 & 1.00 & 1.00 & 1.00 & 1.00 & 1.00 \\
   \hline
\end{tabular}
\end{center}
\caption{Numbers represent probability of rejection under the true model, with rejection region calculated
under $H_0$. $\alpha$-divergences, for $\Delta_n=0.001$}
\label{tabcir1a}
\end{table}

\begin{table}
\begin{center}
\begin{tabular}{lrrrrrr}
  \hline
model ($\alpha$, $n$) & $\lambda=-0.99$ & $\lambda=-1.20$ & $\lambda=-1.50$ & $\lambda=-1.75$ & $\lambda=-2.00$ & $\lambda=-2.50$ \\
  \hline
CIR$_0$ (0.01, 50) & 0.00 & 0.00 & 0.00 & 0.02 & 0.05 & 0.13 \\
  CIR$_1$ (0.01, 50) & 0.00 & 0.01 & 0.28 & 0.55 & 0.71 & 0.87 \\
  CIR$_2$ (0.01, 50) & 0.00 & 0.81 & 1.00 & 1.00 & 1.00 & 1.00 \\
   &  &  &  &  &  &  \\
  CIR$_0$ (0.05, 50) & 0.00 & 0.00 & 0.01 & 0.04 & 0.09 & 0.19 \\
  CIR$_1$ (0.05, 50) & 0.00 & 0.09 & 0.48 & 0.70 & 0.81 & 0.92 \\
  CIR$_2$ (0.05, 50) & 0.00 & 0.97 & 1.00 & 1.00 & 1.00 & 1.00 \\
   &  &  &  &  &  &  \\
   &  &  &  &  &  &  \\
  CIR$_0$ (0.01, 100) & 0.00 & 0.00 & 0.00 & 0.02 & 0.05 & 0.14 \\
  CIR$_1$ (0.01, 100) & 0.00 & 0.21 & 0.83 & 0.95 & 0.98 & 0.99 \\
  CIR$_2$ (0.01, 100) & 0.00 & 1.00 & 1.00 & 1.00 & 1.00 & 1.00 \\
   &  &  &  &  &  &  \\
  CIR$_0$ (0.05, 100) & 0.00 & 0.00 & 0.01 & 0.05 & 0.09 & 0.20 \\
  CIR$_1$ (0.05, 100) & 0.00 & 0.53 & 0.93 & 0.98 & 0.99 & 1.00 \\
  CIR$_2$ (0.05, 100) & 0.24 & 1.00 & 1.00 & 1.00 & 1.00 & 1.00 \\
   &  &  &  &  &  &  \\
   &  &  &  &  &  &  \\
  CIR$_0$ (0.01, 500) & 0.00 & 0.00 & 0.00 & 0.02 & 0.04 & 0.10 \\
  CIR$_1$ (0.01, 500) & 0.00 & 1.00 & 1.00 & 1.00 & 1.00 & 1.00 \\
  CIR$_2$ (0.01, 500) & 1.00 & 1.00 & 1.00 & 1.00 & 1.00 & 1.00 \\
   &  &  &  &  &  &  \\
  CIR$_0$ (0.05, 500) & 0.00 & 0.00 & 0.01 & 0.04 & 0.07 & 0.15 \\
  CIR$_1$ (0.05, 500) & 0.95 & 1.00 & 1.00 & 1.00 & 1.00 & 1.00 \\
  CIR$_2$ (0.05, 500) & 1.00 & 1.00 & 1.00 & 1.00 & 1.00 & 1.00 \\
   \hline
\end{tabular}
\end{center}
\caption{Numbers represent probability of rejection under the true
model, with rejection region calculated under $H_0$.
Power-divergences  for $\Delta_n=0.001$} \label{tabcir1p}
\end{table}

\begin{table}
\begin{center}
\begin{tabular}{lrrrrrr}
  \hline
model ($\alpha$, $n$) & $\alpha=-0.99$ & $\alpha=-0.90$ & $\alpha=-0.75$ & $\alpha=-0.50$ & $\alpha=-0.25$ & $\alpha=-0.10$ \\
  \hline
CIR$_0$ (0.01, 50) & 0.01 & 0.14 & 0.54 & 0.77 & 0.85 & 0.87 \\
  CIR$_1$ (0.01, 50) & 0.80 & 0.98 & 1.00 & 1.00 & 1.00 & 1.00 \\
  CIR$_2$ (0.01, 50) & 1.00 & 1.00 & 1.00 & 1.00 & 1.00 & 1.00 \\
   &  &  &  &  &  &  \\
  CIR$_0$ (0.05, 50) & 0.05 & 0.16 & 0.54 & 0.77 & 0.85 & 0.87 \\
  CIR$_1$ (0.05, 50) & 0.94 & 0.98 & 1.00 & 1.00 & 1.00 & 1.00 \\
  CIR$_2$ (0.05, 50) & 1.00 & 1.00 & 1.00 & 1.00 & 1.00 & 1.00 \\
   &  &  &  &  &  &  \\
   &  &  &  &  &  &  \\
  CIR$_0$ (0.01, 100) & 0.01 & 0.13 & 0.49 & 0.71 & 0.79 & 0.82 \\
  CIR$_1$ (0.01, 100) & 0.99 & 1.00 & 1.00 & 1.00 & 1.00 & 1.00 \\
  CIR$_2$ (0.01, 100) & 1.00 & 1.00 & 1.00 & 1.00 & 1.00 & 1.00 \\
   &  &  &  &  &  &  \\
  CIR$_0$ (0.05, 100) & 0.04 & 0.15 & 0.49 & 0.71 & 0.79 & 0.82 \\
  CIR$_1$ (0.05, 100) & 1.00 & 1.00 & 1.00 & 1.00 & 1.00 & 1.00 \\
  CIR$_2$ (0.05, 100) & 1.00 & 1.00 & 1.00 & 1.00 & 1.00 & 1.00 \\
   &  &  &  &  &  &  \\
   &  &  &  &  &  &  \\
  CIR$_0$ (0.01, 500) & 0.00 & 0.06 & 0.28 & 0.54 & 0.69 & 0.74 \\
  CIR$_1$ (0.01, 500) & 1.00 & 1.00 & 1.00 & 1.00 & 1.00 & 1.00 \\
  CIR$_2$ (0.01, 500) & 1.00 & 1.00 & 1.00 & 1.00 & 1.00 & 1.00 \\
   &  &  &  &  &  &  \\
  CIR$_0$ (0.05, 500) & 0.02 & 0.08 & 0.28 & 0.54 & 0.69 & 0.74 \\
  CIR$_1$ (0.05, 500) & 1.00 & 1.00 & 1.00 & 1.00 & 1.00 & 1.00 \\
  CIR$_2$ (0.05, 500) & 1.00 & 1.00 & 1.00 & 1.00 & 1.00 & 1.00 \\
   \hline
\end{tabular}
\end{center}
\caption{Numbers represent probability of rejection under the true model, with rejection region calculated
under $H_0$. $\alpha$-divergences, for $\Delta_n=0.1$}
\label{tabcir3a}
\end{table}

\begin{table}
\begin{center}
\begin{tabular}{lrrrrrr}
  \hline
model ($\alpha$, $n$) & $\lambda=-0.99$ & $\lambda=-1.20$ & $\lambda=-1.50$ & $\lambda=-1.75$ & $\lambda=-2.00$ & $\lambda=-2.50$ \\
  \hline
CIR$_0$ (0.01, 50) & 0.00 & 0.00 & 0.00 & 0.01 & 0.02 & 0.06 \\
  CIR$_1$ (0.01, 50) & 0.00 & 0.06 & 0.52 & 0.75 & 0.86 & 0.94 \\
  CIR$_2$ (0.01, 50) & 0.00 & 0.99 & 1.00 & 1.00 & 1.00 & 1.00 \\
   &  &  &  &  &  &  \\
  CIR$_0$ (0.05, 50) & 0.00 & 0.00 & 0.00 & 0.02 & 0.04 & 0.09 \\
  CIR$_1$ (0.05, 50) & 0.00 & 0.23 & 0.70 & 0.85 & 0.92 & 0.96 \\
  CIR$_2$ (0.05, 50) & 0.06 & 1.00 & 1.00 & 1.00 & 1.00 & 1.00 \\
   &  &  &  &  &  &  \\
   &  &  &  &  &  &  \\
  CIR$_0$ (0.01, 100) & 0.00 & 0.00 & 0.00 & 0.01 & 0.02 & 0.05 \\
  CIR$_1$ (0.01, 100) & 0.00 & 0.56 & 0.96 & 0.99 & 1.00 & 1.00 \\
  CIR$_2$ (0.01, 100) & 0.00 & 1.00 & 1.00 & 1.00 & 1.00 & 1.00 \\
   &  &  &  &  &  &  \\
  CIR$_0$ (0.05, 100) & 0.00 & 0.00 & 0.00 & 0.02 & 0.03 & 0.08 \\
  CIR$_1$ (0.05, 100) & 0.00 & 0.83 & 0.99 & 1.00 & 1.00 & 1.00 \\
  CIR$_2$ (0.05, 100) & 0.97 & 1.00 & 1.00 & 1.00 & 1.00 & 1.00 \\
   &  &  &  &  &  &  \\
   &  &  &  &  &  &  \\
  CIR$_0$ (0.01, 500) & 0.00 & 0.00 & 0.00 & 0.00 & 0.01 & 0.02 \\
  CIR$_1$ (0.01, 500) & 0.00 & 1.00 & 1.00 & 1.00 & 1.00 & 1.00 \\
  CIR$_2$ (0.01, 500) & 1.00 & 1.00 & 1.00 & 1.00 & 1.00 & 1.00 \\
   &  &  &  &  &  &  \\
  CIR$_0$ (0.05, 500) & 0.00 & 0.00 & 0.00 & 0.01 & 0.02 & 0.04 \\
  CIR$_1$ (0.05, 500) & 1.00 & 1.00 & 1.00 & 1.00 & 1.00 & 1.00 \\
  CIR$_2$ (0.05, 500) & 1.00 & 1.00 & 1.00 & 1.00 & 1.00 & 1.00 \\
   \hline
\end{tabular}
\end{center}
\caption{Numbers represent probability of rejection under the true
model, with rejection region calculated under $H_0$.
Power-divergences  for $\Delta_n=0.1$} \label{tabcir3p}
\end{table}
\end{landscape}

\end{document}